\documentclass[11pt]{article}
\usepackage{amssymb,amsfonts,amsmath,amsthm, psfrag,eepic,amscd}
\newtheorem{theo}{Theorem}

\theoremstyle{definition}

\newtheorem{exam}[theo]{Example}

\makeatletter \@addtoreset{equation}{} \@addtoreset{theo}{section}
\makeatother

\begin{document}
\title{\bf Generalizations of Han's Hook Length Identities}
\author{Laura L.M. Yang\\
\small Department of Mathematical and Statistical Sciences\\
\small University of Alberta, Edmonton, Alberta, Canada T6G 2G1\\
\small \texttt{yanglm@hotmail.com} }
\date{ } 
\maketitle
\begin{abstract}
\noindent Han recently discovered new hook length identities for
binary trees.  In this paper, we extend Han's identities to binomial
families of trees. Moreover, we present a bijective proof of one of
the identities for the family of ordered trees.
\end{abstract}

\section{Introduction}
The {\it hook length} of a vertex $v$ of a rooted tree $T$ is the
number $h_v$ of descendants of $v$ in $T$ (including $v$ itself).
Several identities involving this parameter have been discovered,
especially since the appearance of Postnikov's identity in 2004
\cite{Postnikov_y2004}; see, e.g.,
\cite{Chen_y2008,Du_y2005,Gessel_y2005,Moon_y2007,Seo_y2005} and the
references contained therein. Han \cite{Han_1,Han_2} recently found
two more such identities, namely,
\begin{eqnarray}
\sum_{T}\prod_{v\in T}\frac{1}{h_v2^{h_v-1}} &=&\frac{1}{n!}
\end{eqnarray}
and
\begin{eqnarray}
\sum_{T}\prod_{v\in T}\frac{(z+h_v)^{h_v-1}}{h_v(2z+h_v-1)^{h_v-2}}
&=&\frac{2^nz}{n!}(z+n)^{n-1},
\end{eqnarray}
where each sum is over the (incomplete) binary trees $T$ with $n$
vertices (in which each vertex has at most one left-child and at
most one right-child). Our main object here is to extend Han's
identities to more general {\it binomial} families of trees. The
definition of these families and our main results will be stated in
Section 2. The proofs will be given in Section 3. Finally, in
Section 4, we give a bijective proof of one of the identities for
the family of ordered trees.
\section{Definitions and Main Results}
We recall that {\it ordered} trees are (finite) rooted trees with an
ordering specified for the children of each vertex (see, e.g., Knuth
\cite[p. 306]{Knuth_y1973}). Let $s$ and $m$ be given constants such
that $sm>0$ and $m$ is a positive integer if $s>0$. And let $d_v$
denote the number of children of vertex $v$ in any given rooted tree
$T$. If each ordered tree $T$ is assigned the {\it weight}
$$
w(T)=\prod_{v\in T} \binom{m}{d_v}s^{d_v},
$$
then the resulting family of weighted ordered trees is called a {\it
binomial} family or an $(s,m)$-family $F$. Let $F_n$ denote the
subset of binomial trees that have $n$ vertices and let
$y_n=\sum_{T\in F_n}w(T)$ denote the (weighted) number of trees in
$F_n$. It follows readily from these definitions that the generating
function $y=y(x)=\sum_1^\infty y_nx^n$ satisfies the relation
$$y=x(1+sy)^m.$$
For additional remarks on these families, especially in the context
of simply generated families may be sound in \cite{Moon_y2007}.
Notice, for example, that the binomial families include the
incomplete k-ary and the ordered trees; but they do not include the
complete binary trees in which every vertex has zero or two
children. We now state our main results.

\begin{theo}
Let $F_n$ denote the subset of the $(s,m)$-family of binomial trees
that have $n$ vertices. Then
\begin{eqnarray}\label{Identity-1}
\sum_{T\in F_n}w(T)\prod_{v\in T}
\frac{1}{h_vm^{h_v-1}}&=&\frac{s^{n-1}}{n!}
\end{eqnarray}
and
\begin{eqnarray}\label{Identity-2}
\sum_{T\in F_n}w(T)\prod_{v\in T}
\frac{(z+h_v)^{h_v-1}}{h_v(mz+h_v-1)^{h_v-2}}
&=&\frac{s^{n-1}m^nz}{n!}(z+n)^{n-1},
\end{eqnarray}
for $n=1,2,\ldots$.
\end{theo}

\section{Proof of Theorem}
Let $p_n$ and $q_n$ denote the lefthand sides of identities
\eqref{Identity-1} and \eqref{Identity-2} for $n=1,2,\ldots$. The
proof will be by induction on $n$. It is easy to check that $p_1=1$
and $q_1=mz$ so \eqref{Identity-1} and \eqref{Identity-2} hold when
$n=1$. Any non-trivial binomial tree $T$ with $n$ vertices in which
the root has $d$ children may be constructed from an ordered
collection of d smaller binomial trees with $n-1$ vertices
altogether by attaching a new (root) vertex to the roots of the $d$
smaller trees and then introducing the appropriate weight factors.
It follows readily from this observation and the definition of
$p_n$, that if $n>1$ then
\begin{eqnarray}\label{p_n}
p_n=\frac{1}{nm^{n-1}}\sum_{d\geq 1}\binom{m}{d}s^d\sum
p_{j_1}\ldots p_{j_d}
\end{eqnarray}
where the inner sum is over all compositions $(j)=(j_1,\ldots,j_d)$
of $n-1$ into $d$ positive integers. If we apply the induction
hypothesis that $p_j=s_{j-1}/j!$ for $j<n$, simplify, and rewrite
the righthand side of relation \eqref{p_n} in terms of generating
functions, we find that
\begin{eqnarray*}
p_n &=&\frac{s^{n-1}}{nm^{n-1}}\sum_{d\geq 1}\binom{m}{d}
[x^{n-1}]{(e^x-1)^d}\\
&=&\frac{s^{n-1}}{nm^{n-1}}[x^{n-1}](e^{xm}-1) = \frac{s^{n-1}}{n!}.
\end{eqnarray*}
This suffices to prove identity \eqref{Identity-1}.

Before proceeding to the proof of identity \eqref{Identity-2} we
recall that if $u=u(x)$ is a power series such that $u=e^{xu}$, then
it follows readily from Lagrange's inversion formula that
\begin{eqnarray}\label{u}
u^z=1+\sum_{n\geq 1}\frac{z(z+n)^{n-1}}{n!}x^n.
\end{eqnarray}
for any $z$.

We now consider identify \eqref{Identity-2} for the quantity $q_n$.
In this case the reasoning that led to relation \eqref{p_n} leads to
the conclusion that if $n>1$, then
\begin{eqnarray}\label{q_n}
q_n=\frac{(z+n)^{n-1}}{n(mz+n-1)^{n-2}}\sum_{d\geq 1}
\binom{m}{d}s^d \sum q_{j_1}\ldots q_{j_d}
\end{eqnarray}
where the inner sum is over the same compositions $(j)$ as before.
If we apply the induction hypothesis that
$q_j=s^{j-1}m^jz(z+j)^{j-1}/j!$ for $j<n$, simplify, rewrite the
righthand side of relation \eqref{q_n} in terms of generating
functions, and appeal to relation \eqref{u}, we find that
\begin{eqnarray*}
q_n
&=&\frac{(sm(z+n))^{n-1}}{n(mz+n-1)^{n-2}}\sum_{d\geq1}\binom{m}{d}[x^{n-1}]{(u^z-1)^d}\\
&=&\frac{(sm(z+n))^{n-1}}{n(mz+n-1)^{n-2}}[x^{n-1}](u^{zm}-1) =
\frac{s^{n-1}m^n z(z+n)^{n-1}}{n!}.
\end{eqnarray*}
This suffices to complete the proof of the theorem.

\begin{exam}
The five ordered trees with $n=4$ vertices are illustrated in Figure
\ref{eg-ordered}.
\begin{figure}[h]
\begin{center}
\begin{picture}(370,60)
\setlength{\unitlength}{0.8mm}

\linethickness{0.4pt} \put(5,15){\circle*{1}}
\put(5,15){\line(-1,-1){10}} \put(-5,5){\circle*{1}}
\put(5,15){\line(0,-1){10}} \put(5,5){\circle*{1}}
\put(5,15){\line(1,-1){10}} \put(15,5){\circle*{1}}
\put(4,20){$T_1$}

\put(45,15){\circle*{1}} \put(45,15){\line(-1,-1){10}}
\put(35,5){\circle*{1}} \put(45,15){\line(1,-1){10}}
\put(55,5){\circle*{1}} \put(35,5){\line(0,-1){10}}
\put(35,-5){\circle*{1}} \put(44,20){$T_2$}

\put(85,15){\circle*{1}} \put(85,15){\line(-1,-1){10}}
\put(75,5){\circle*{1}} \put(85,15){\line(1,-1){10}}
\put(95,5){\circle*{1}} \put(95,5){\line(0,-1){10}}
\put(95,-5){\circle*{1}} \put(84,20){$T_3$}

\put(125,15){\circle*{1}} \put(125,15){\line(0,-1){10}}
\put(125,5){\circle*{1}} \put(125,5){\line(-1,-1){10}}
\put(115,-5){\circle*{1}} \put(125,5){\line(1,-1){10}}
\put(135,-5){\circle*{1}} \put(124,20){$T_4$}

\put(165,15){\circle*{1}} \put(165,15){\line(0,-1){20}}
\put(165,15){\circle*{1}} \put(165,9){\circle*{1}}
\put(165,3){\circle*{1}} \put(165,-5){\circle*{1}}
\put(164,20){$T_5$}
\end{picture}
\end{center}
\caption{Ordered trees with $4$ vertices.} \label{eg-ordered}
\end{figure}
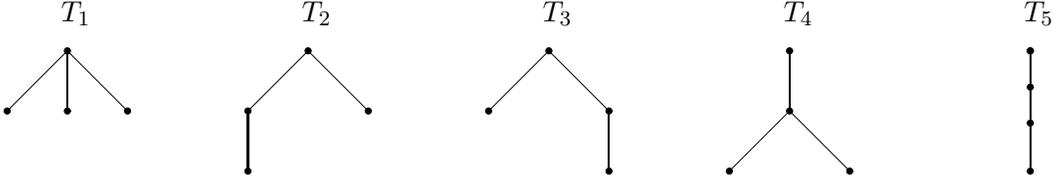
If $F$ is the $(1,k)$-family (of incomplete $k$-ary trees), then it
follows from the Theorem that
\begin{eqnarray}
\sum_{T}w(T)\prod_{v\in T}
\frac{1}{h_vk^{h_v-1}}&=&\frac{1}{n!}\label{(1,k)-1}
\end{eqnarray}
and
\begin{eqnarray}
\sum_{T}w(T)\prod_{v\in T}
\frac{(z+h_v)^{h_v-1}}{h_v(kz+h_v-1)^{h_v-2}}
&=&\frac{k^nz}{n!}(z+n)^{n-1},\label{(1,k)-2}
\end{eqnarray}
where the sums, here and elsewhere, are over the trees $T$ in $F_n$.
In this case, the weights of $T_1, T_2, T_3, T_4$ and $T_5$ are
$\binom{k}{3}, \binom{k}{2}\binom{k}{1}, \binom{k}{2}\binom{k}{1},
\binom{k}{1}\binom{k}{2}$ and $\left(\binom{k}{1}\right)^3$,
respectively. Hence,
\begin{eqnarray*}
p_4=\frac{\binom{k}{3}}{4\cdot
k^3}+\frac{\binom{k}{2}\binom{k}{1}}{4\cdot k^3 \cdot 2\cdot
k}+\frac{\binom{k}{2}\binom{k}{1}}{4\cdot k^3 \cdot 2\cdot
k}+\frac{\binom{k}{1}\binom{k}{2}}{4\cdot k^3\cdot 3\cdot
k^2}+\frac{\left(\binom{k}{1}\right)^3}{4\cdot k^3 \cdot 3\cdot
k^2\cdot 2\cdot k}=\frac{1}{4!}
\end{eqnarray*}
and
\begin{eqnarray*}
q_4&=&\frac{\binom{k}{3}(z+4)^3}{4(kz+3)^2(kz)^{-3}}
+\frac{\binom{k}{2}\binom{k}{1}(z+4)^3(z+2)}{4(kz+3)^2\cdot
2(kz)^{-2}}
+\frac{\binom{k}{2}\binom{k}{1}(z+4)^3(z+2)}{4(kz+3)^2\cdot 2(kz)^{-2}}\\
&&+\frac{\binom{k}{1}\binom{k}{2}(z+4)^3(z+3)^2}{4(kz+3)^2\cdot
3(kz+2)(kz)^{-2}}
+\frac{\left(\binom{k}{1}\right)^3(z+4)^3(z+3)^2(z+2)}{4(kz+3)^2\cdot3(kz+2)\cdot2(kz)^{-1}}
=\frac{k^4z(z+4)^3}{4!}.
\end{eqnarray*}
Notice that \eqref{(1,k)-1} and \eqref{(1,k)-2} reduce to Han's
identities when $k=2$.

If $F$ is a $(-1,-k)$-family, then
\begin{eqnarray}
\sum_{T}w(T)\prod_{v\in T}
\frac{1}{h_v(-k)^{h_v-1}}&=&\frac{(-1)^{n-1}}{n!}\label{(-1,-k)-family}
\end{eqnarray}
and
\begin{eqnarray}
\sum_{T}w(T)\prod_{v\in T}
\frac{(z+h_v)^{h_v-1}}{h_v(h_v-kz-1)^{h_v-2}}
&=&\frac{-k^nz}{n!}(z+n)^{n-1}.
\end{eqnarray}
In particular, if $F$ is the $(-1,-1)$-family, i.e., the family of
ordered trees, then
\begin{eqnarray}
\sum_{T}\prod_{v\in T}
\frac{1}{h_v(-1)^{h_v-1}}&=&\frac{(-1)^{n-1}}{n!}\label{ordered}
\end{eqnarray}
and
\begin{eqnarray}
\sum_{T}\prod_{v\in T} \frac{(z+h_v)^{h_v-1}}{h_v(h_v-z-1)^{h_v-2}}
&=&\frac{-z}{n!}(z+n)^{n-1},
\end{eqnarray}
where have omitted the weight factors here since they all equal one.
In this case,
\begin{eqnarray*} p_4&=&\frac{1}{4\cdot
(-1)^3}+\frac{1}{4\cdot (-1)^3 \cdot 2\cdot
(-1)}+\frac{1}{4\cdot (-1)^3 \cdot 2\cdot (-1)}\\
&&+\frac{1}{4\cdot (-1)^3\cdot 3\cdot (-1)^2}+\frac{1}{4\cdot (-1)^3
\cdot 3\cdot (-1)^2\cdot 2\cdot (-1)}=-\frac{1}{4!}
\end{eqnarray*}
and
\begin{eqnarray*}
q_4&=&\frac{(z+4)^3}{4(3-z)^2(-z)^{-3}}
+\frac{(z+4)^3(z+2)}{4(3-z)^2\cdot 2(-z)^{-2}}
+\frac{(z+4)^3(z+2)}{4(3-z)^2\cdot 2(-z)^{-2}}\\
&&+\frac{(z+4)^3(z+3)^2}{4(3-z)^2\cdot 3(2-z)(-z)^{-2}}
+\frac{(z+4)^3(z+3)^2(z+2)}{4(3-z)^2\cdot3(2-z)\cdot2(-z)^{-1}}
=-\frac{z(z+4)^3}{4!}.
\end{eqnarray*}

If $F$ is a $(1/m,m)$-family, then
\begin{eqnarray}
\sum_{T}w(T)\prod_{v\in T}
\frac{1}{h_vm^{h_v-1}}&=&\frac{1}{m^{n-1}n!}
\end{eqnarray}
and
\begin{eqnarray}
\sum_{T}w(T)\prod_{v\in T}
\frac{(z+h_v)^{h_v-1}}{h_v(mz+h_v-1)^{h_v-2}}
&=&\frac{mz}{n!}(z+n)^{n-1} \label{rooted tree}.
\end{eqnarray}
If we let $z=1/m$ in \eqref{rooted tree} and take the limit as $m$
tends to infinity, we obtain the identity
\begin{eqnarray}
\sum_{T} \prod_{v\in T}\frac{1}{d_v!}=\frac{n^{n-1}}{n!},
\end{eqnarray}
where the sum is over all ordered trees $T$ with $n$ vertices. This
relation, which expresses the number $n^{n-1}$ of rooted labelled
trees with $n$ vertices as a sum over the ordered trees with $n$
vertices, with suitable weights taken into account, is equivalent to
a relation given by Mohanty \cite[p. 163]{Mohanty_y1979}.
\end{exam}

\section{An Involution on Increasing Ordered Trees}
We conclude by giving a sign-reversing involution that establishes
an alternate form of identity \eqref{ordered}, namely,
\begin{eqnarray}
\sum_{T}n!\prod_{v\in T} \frac{1}{h_v(-1)^{h_v}}&=&-1,
\end{eqnarray}
where the sum is over all ordered trees with $n$ vertices.

It is well known that $n!/\prod_{v\in T}h_v$ counts the number of
ways to label the vertices of $T$ with $\{1,2,\ldots, n\}$ such that
the label of each vertex is less than that of its descendants
\cite[p.67, exer. 20]{Knuth_y1997}. Such a labelled tree is called
{\it increasing}. We define the sign of a tree $T$ to be
$\prod_{v\in T} (-1)^{h_v}$.

An increasing ordered tree $T$ with $n$ vertices is {\it proper} if
the root of $T$ has $n-1$ children and their labels are increasing
from left to right. It is easy to check that the sign of any proper
tree is $-1$.

The involution is based on the non-proper increasing ordered trees.
For any leaf $v$ of a non-proper increasing ordered tree, suppose
$v$ is the $i$-th child of $u$ and $w$ is the $i+1$-th child of $u$
if it exists. We say that $v$ is {\it illegal} if $v$ is the
rightmost child of $u$ and the subtree rooted at $u$ is proper or
$v$ is bigger than any vertex of the subtree rooted at $w$ and the
subtree rooted at $w$ is proper.

Now the involution can be described as follows: Given any
non-proper increasing ordered tree $T$ let $v$ be the first
illegal leaf encountered when traversing the tree $T$ in preorder.
We now have two cases: (1) $v$ is bigger than any vertex of the
subtree rooted at $w$ and the subtree rooted at $w$ is proper; (2)
$v$ is the rightmost child of $u$ and the subtree rooted at $u$ is
proper. In this case, $u$ is not the root of $T$.

For case (1), let $u$ be the parent of $v$.  We cut off the edge
between $u$ and $v$, and move $v$ as the rightmost child of $w$. Let
$T'$ be the resulting tree. Note that in the search process for
$T'$, the leaf $v$ is still the first encountered illegal leaf.

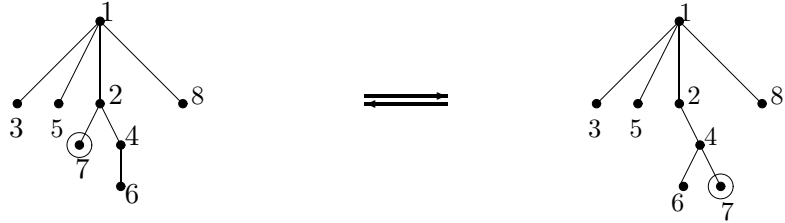
\begin{figure}[h,t]
\begin{center}
\begin{picture}(340,80)
\setlength{\unitlength}{1.1mm}

\linethickness{0.4pt} \put(5,15){\circle*{1}}
\put(15,25){\circle*{1}} \put(15,25){\line(-1,-1){10}}
\put(4,11){3} \put(15,25){1}

\put(15,25){\line(0,-1){10}} \put(15,15){\circle*{1}}
\put(16,15){2}

\put(15,15){\line(-1,-2){2.5}} \put(15,15){\line(1,-2){2.5}}
\put(17.5,10){\line(0,-1){5}} \put(12.5,10){\circle{3}}
\put(12.5,10){\circle*{1}} \put(17.5,10){\circle*{1}}
\put(17.5,5){\circle*{1}} \put(12,6){7}
\put(18,10){4}\put(18,3){6}
 \small
 \put(15,25){\line(-1,-2){5}} \put(10,15){\circle*{1}}
\put(9,11){5}

\put(15,25){\line(1,-1){10}} \put(25,15){\circle*{1}}
\put(26,15){8}

\linethickness{1pt} \put(47,16){\vector(1,0){10}}
\put(57,15){\vector(-1,0){10}}
\linethickness{0.4pt} \put(75,15){\circle*{1}}
\put(85,25){\circle*{1}} \put(85,25){\line(-1,-1){10}}
\put(74,11){3} \put(85,25){1}

\put(85,25){\line(0,-1){10}} \put(85,15){\circle*{1}}
\put(86,15){2}

\put(85,15){\line(1,-2){2.5}} \put(87.5,10){\line(-1,-2){2.5}}
\put(87.5,10){\line(1,-2){2.5}} \put(87.5,10){\circle*{1}}
\put(85.5,5){\circle*{1}} \put(90,5){\circle*{1}}
\put(90,5){\circle{3}} \put(88,10){4} \put(84,2){6} \put(90,1){7}

\put(85,25){\line(-1,-2){5}} \put(80,15){\circle*{1}}
\put(79,11){5}

\put(85,25){\line(1,-1){10}} \put(95,15){\circle*{1}}
\put(96,15){8}

\end{picture}
\end{center}
\caption{An involution on increasing ordered trees.}
\label{eg-involution}
\end{figure}

For case (2), we may reverse the construction for case (1). Hence
we obtain a sign-reversing involution. Figure \ref{eg-involution}
illuminates this involution on increasing ordered trees.


\end{document}